\documentclass[12pt,a4paper]{article}
\usepackage{amssymb}

\usepackage{amssymb}
\usepackage{hyperref}

\usepackage{graphicx,amssymb,amsfonts,epsfig,a4,amsmath,amsthm,xypic}
\usepackage[latin1]{inputenc}

\usepackage{tikz-cd}

\newtheorem{Thm}{Theorem}[section]
\newtheorem{Prop}[Thm]{Proposition}
\newtheorem{Lem}[Thm]{Lemma}
\newtheorem{Cor}[Thm]{Corollary}

\newtheorem{Rem}[Thm]{Remark}
\newtheorem{Ex}{Example}

\newcommand{\Z}{\mathbb{Z}}

\newcommand{\N}{\mathbf{N}}
\newcommand{\R}{\mathbb{R}}
\newcommand{\C}{\mathbb{C}}
\newcommand{\Q}{\mathbb{Q}}

\newcommand{\F}{\mathbf{F}}

\newcommand{\Ker}{\text{Ker}}

\newcommand{\bpr}{\noindent \textbf{Proof}: }

\newcommand{\epr}{~$\blacksquare$}
\newcommand{\eprsk}{~$\blacksquare$\medskip}

\title{Maximal Haagerup subgroups in $\Z^{n+1}\rtimes_{\rho_n} GL_2(\Z)$}
\author{Alain VALETTE }

\begin{document}

\maketitle

\begin{abstract} For $n\geq 1$, let $\rho_n$ denote the standard action of $GL_2(\Z)$ on the space $P_n(\Z)\simeq\Z^{n+1}$ of homogeneous polynomials of degree $n$ in two variables, with integer coefficients. For $G$ a non-amenable subgroup of $GL_2(\Z)$, we describe the maximal Haagerup subgroups of the semi-direct product $\Z^{n+1}\rtimes_{\rho_n} G$, extending the classification of Jiang-Skalski \cite{JiSk} of the maximal Haagerup subgroups in $\Z^2\rtimes SL_2(\Z)$. We prove that, for $n$ odd, the group $P_n(\Z)\rtimes SL_2(\Z)$ admits infinitely many pairwise non-conjugate maximal Haagerup subgroups which are free groups; and that, for $n$ even, the group $P_n(\Z)\rtimes GL_2(\Z)$ admits infinitely many pairwise non-conjugate maximal Haagerup subgroups which are isomorphic to $SL_2(\Z)$.\end{abstract}

\section{Introduction}

For discrete countable groups, the Haagerup property is a weak form of amenability that proved to be useful in many questions in analytical group theory, ranging from K-theory to dynamical systems (see \cite{CCJJV}). It is not difficult to see that, in a countable group, every Haagerup subgroup is contained in a maximal one (see Proposition 1.3 in \cite{JiSk} or lemma \ref{maxHaag} below). This raises the question, given a group $G$, of describing the maximal Haagerup subgroups of $G$.

The study of maximal Haagerup subgroups was initiated by Y. Jiang and A. Skalski \cite{JiSk}, and we refer to this paper for many interesting and intriguing examples. We mention here Theorem 2.12 in \cite{JiSk}, where the authors classify maximal Haagerup subgroups of the semi-direct product $\Z^2\rtimes SL_2(\Z)$. This example is especially interesting in view of a result of Burger (Example 2 following Proposition 7 in \cite{Burger}): if $G$ is a non-amenable subgroup of $SL_2(\Z)$, then the pair $(\Z^2\rtimes G,\Z^2)$ has the relative property (T); in particular $\Z^2\rtimes G$ is not Haagerup, in spite of the fact that $\Z^2$ and $G$ are both Haagerup.

\begin{Thm}[Theorem 2.12 in \cite{JiSk}]\label{JiaSka} Let $H$ be a maximal Haagerup subgroup of $\Z^2\rtimes SL_2(\Z)$. Then there is a dichotomy:
\begin{enumerate}
\item either $H=\Z^2\rtimes C$, where $C$ is a maximal amenable\footnote{We classify those subgroups in Proposition \ref{MaxAmen}.}  subgroup of $SL_2(\Z)$ ;
\item or $H\cap\Z^2$ is trivial; then $H$ is not amenable. If $K$ denotes the image of $H$ under the quotient map $\Z^2\rtimes SL_2(\Z)\rightarrow SL_2(\Z)$ (so that $K$ is isomorphic to $H$), then $H=\{(b(g),g):g\in K\}$ where $b:K\rightarrow \Z^2$ is a 1-cocycle that cannot be extended to a larger subgroup of $SL_2(\Z)$.
\end{enumerate}
\end{Thm}

\begin{Rem} Denote by $L(G)$ the group von Neumann algebra of the group $G$. In Theorem 3.1 of \cite{JiSk}, Jiang and Skalski are able to prove the stronger result that, if $C$ is a maximal amenable subgroup of $SL_2(\Z)$, such that $\Z^2\rtimes C$ has infinite conjugacy classes, then $L(\Z^2\rtimes C)$ is a maximal Haagerup von Neumann subalgebra of $L(\Z^2\rtimes SL_2(\Z))$, where the Haagerup property for finite von Neumann algebras was defined in \cite{Jol}. Subsequently Y. Jiang (Corollary 4.3 in \cite{Jiang}) showed that $L(SL_2(\Z))$ is a maximal Haagerup subalgebra in $L(\Z^2\rtimes SL_2(\Z))$. It was pointed out to us by A. Skalski, that there is no known example of a maximal Haagerup subgroup $H$ in a group $G$, such that $L(H)$ is not maximal Haagerup in $L(G)$.
\end{Rem}

We now come to the present paper. Fix $n\geq 1$ and, for $A=\Z,\R,\C$, denote by $P_n(A)$ the set of polynomials in two variables $X,Y$, with coefficients in $A$, which are homogeneous of degree $n$, so that $P_n(A)\simeq A^{n+1}$. It is a classical fact that $GL_2(\R)$ admits an irreducible representation $\rho_n$ on $P_n(\R)$ given as follows: for $P\in P_n(\R)$ and $A=\left(\begin{array}{cc}a_{11} & a_{12} \\a_{21} & a_{22}\end{array}\right)\in GL_2(\R)$, set:
$$(\rho_n(A)(P))(X,Y)= P((X,Y)\cdot A)=P(a_{11}X+a_{21}Y,a_{12}X+a_{22}Y).$$
Since $\rho_n(GL_2(\Z))$ leaves $P_n(\Z)$ invariant, we may form the semi-direct product
$$G_n=:P_n(\Z)\rtimes_{\rho_n} GL_2(\Z)$$
(observe $G_1=\Z^2\rtimes GL_2(\Z)$, that contains $\Z^2\rtimes SL_2(\Z)$ as a subgroup of index 2). Let $G$ be a non-amenable subgroup of $GL_2(\Z)$: our goal is to classify the maximal Haagerup subgroups of $P_n(\Z)\rtimes_{\rho_n} G\subset G_n$. So here is our first main result, extending Theorem \ref{JiaSka}:

\begin{Thm}\label{main} Fix $n\geq 1$ and a non-amenable subgroup $G$ of $GL_2(\Z)$. Let $H$ be a maximal Haagerup subgroup of $P_n(\Z)\rtimes_{\rho_n} G$. Then there is a dichotomy.
\begin{enumerate} 
\item Either $H$ is amenable, in which case $H=P_n(\Z)\rtimes_{\rho_n}C$, with $C$ maximal amenable in $G$;
\item or $H$ is non amenable and there exists a subgroup $K\subset G$ isomorphic to $H$, and a 1-cocycle $b\in Z^1(K,P_n(\Z))$ such that 
$$H=\{(b(k),k):k\in K\}$$
and $b$ cannot be extended to a larger subgroup of $G$. (In particular, if $b$ is a 1-coboundary, then $K=G$.)
\end{enumerate}
Conversely, any subgroup of $P_n(\Z)\rtimes_{\rho_n} G$ of one of the above 2 forms, defines a maximal Haagerup subgroup in $P_n(\Z)\rtimes_{\rho_n} G$.
\end{Thm}

Even if the conclusion looks very similar to Theorem \ref{JiaSka}, we emphasize that we had to come up with a totally different argument to show that the intersection $H\cap P_n(\Z)$ is either $P_n(\Z)$ or trivial. 

Theorems \ref{JiaSka} and \ref{main} raise an interesting - and somewhat unusual - question in the 1-cohomology of groups: to describe maximal Haagerup subgroups that are non-amenable, we need to describe 1-cocycles that cannot be extended to a larger subgroup. In the case of $\rho_1$, this question is extensively studied in section 2 of \cite{JiSk} (from Lemma 2.14 to Proposition 2.18). 
We attack the question by observing that, for $K$ a subgroup in $SL_2(\Z)$ or $GL_2(\Z)$, a cocycle $b$ on $K$ cannot be extended to a larger subgroup if and only if, for any overgroup $L$ of $K$, the 1-cohomology class of $b$ is not in the image of the restriction map $H^1(L,P_n(\Z))\rightarrow H^1(K,P_n(\Z))$. This suggests to look for subgroups $K$ which are maximal, or close to being maximal, in $SL_2(\Z)$ or $GL_2(\Z)$. Using this approach we prove the following result, which seems to be new even for $n=1$:

\begin{Thm}\label{freesbgp} Assume that $n\geq 1$ is odd. Then there exists a free subgroup $K$ in $SL_2(\Z)$ such that $P_n(\Z)\rtimes SL_2(\Z)$ contains infinitely many maximal Haagerup subgroups of the form $H=\{(b(k),k):k\in K\}$, with $b\in Z^1(K,P_n(\Z))$, which are pairwise not conjugate under $P_n(\Z)$. Moreover the subgroup $K$ may be chosen either with infinite index or with arbitrarily large finite index. 
\end{Thm}

In the even case we get:

\begin{Thm}\label{main2} Assume that $n\geq 2$ is even. Then the semi-direct product $G_n$ contains infinitely many maximal Haagerup subgroups $H$ of the form $H=\{(b(g),g):g\in SL_2(\Z)\}$, where $b\in Z^1(SL_2(\Z),P_n(\Z))$, which moreover are pairwise not conjugate under $P_n(\Z)$.
\end{Thm}

It turns out that the cases of odd and even $n$'s are very different: the reason is that the $GL_2(\Z)$-action on $P_n(\Z)$ factors through $PGL_2(\Z)$ for $n$ even, while it doesn't for $n$ odd.
 


The paper is organized as follows. Section 2 is devoted to prerequisites, with the exception of Proposition \ref{MaxAmen}, describing precisely the maximal amenable subgroups in $SL_2(\Z)$. In section 3 we prove the analogue of Burger's aforementioned result, namely that the pair $(G_n,P_n(\Z))$ has the relative property (T) for $n\geq 1$. Theorem \ref{main} is proved in section 4, while cohomological questions are treated in sections 5 ($n$ odd) and 6 ($n$ even); in particular Theorem \ref{freesbgp} is proved in section 5, and Theorem \ref{main2} in section 6. Theorem \ref{main2} is actually proved by establishing explicit formulae for the ranks \footnote{By the rank of a finitely generated abelian group, we mean the torsion-free rank.} of $H^1(SL_2(\Z),P_n(\Z))$ and $H^1(GL_2(\Z),P_n(\Z))$; those formulae may have their own interest. A combinatorial application of the methods ends the paper.

\bigskip
{\bf Acknowledgements:} Thanks are due to Y. Benoist, Y. Cornulier, L. Hayez, P. Jolissaint, A. Skalski, Y. Stalder and A. Zumbrunnen for useful discussions at various stages of the project. 

\section{Generalities}

\subsection{Haagerup property}

The following results on countable groups will be used freely:
\begin{itemize}
\item The group $GL_2(\Z)$ is Haagerup (see sections 1.2.2 or 1.2.3 in \cite{CCJJV}).
\item An extension of a Haagerup group by an amenable group, is Haagerup; in particular every amenable group is Haagerup (see Proposition 6.1.5 in \cite{CCJJV}).
\end{itemize}

\begin{Lem}\label{maxHaag} Let $G$ be a countable group. Any Haagerup subgroup of $G$ is contained in a maximal Haagerup subgroup.
\end{Lem}

\bpr (compare with Proposition 1.3 in \cite{JiSk}) To apply Zorn's lemma, we must show that being Haagerup is stable by arbitrary increasing unions of subgroups. This follows from the fact that, for countable groups, being Haagerup is a local property, i.e a countable group is Haagerup if and only if every finitely generated subgroup is Haagerup (see Proposition 6.1.1 in \cite{CCJJV}). 
\eprsk

\subsection{1-cohomology of groups}

We recall the following facts from the cohomology of groups. For $G$ a group and $A$ a $G$-module, define the group of 1-cocycles:
$$Z^1(G,A)=\{b:G\rightarrow A:  b(gh)=gb(h)+b(g) \,\mbox{for all}\,g,h\in G \};$$
the group of 1-coboundaries:
$$B^1(G,A)=\{b\in Z^1(G,A): \mbox{there exists}\,a\in A\,\mbox{such that}\,b(g)=ga-a\,\mbox{for all}\,g\in G\};$$
and the first cohomology group:
$$H^1(G,A)=Z^1(G,A)/B^1(G,A).$$
The following is proved e.g. in Proposition 2.3 of \cite{Brown}.

\begin{Prop}\label{1cohom} Let $G$ be a group and let $A$ be a $G$-module. Splittings of the split extension 
$$0\rightarrow A \rightarrow A\rtimes G \rightarrow G \rightarrow 1$$
are given by $i_b:G\rightarrow A\rtimes G:g\mapsto (b(g),g)$ with $b\in Z^1(G,A)$ and are classified up to $A$-conjugacy by the first cohomology group $H^1(G,A)$.
\epr
\end{Prop}

The first part of the following lemma is a variation on Corollaire 7.2, page 39 of \cite{Guich}.

\begin{Lem}\label{central} Let $G$ be a group and $A$ be a $G$-module with cancellation by 2 (i.e. $2x=0\Rightarrow x=0$). 
\begin{enumerate} 
\item Assume $G$ contains a central element $z$ that acts on $A$ by $-1$. Then $H^1(G,A)$ is a vector space over the field with 2 elements. Moreover $b\in B^1(G,A)$ if $b(z)$ belongs to $2A$. 
\item Assume that $G$ contains a central element $c$ of order 2, such that the module action on $A$ factors through $G/\langle c\rangle$. Then the map $H^1(G/\langle c\rangle, A)\rightarrow H^1(G,A)$ (induced by the quotient map $G\rightarrow G/\langle c\rangle$) is an isomorphism. 
\end{enumerate}

\end{Lem}

\bpr \begin{enumerate}
\item We must prove that, for every $b\in Z^1(G,A)$, we have $2b\in B^1(G,A)$. But for $g\in G$ we have, using $zg=gz$ in the 3rd equality:
$$b(z)-b(g)=b(z)+zb(g)=b(zg)=b(gz)=gb(z)+b(g).$$
Re-arranging:
$$2b(g)=(1-g)b(z).$$
Hence $2b\in B^1(G,A)$. If we may write $b(z)=2a$ for some $a\in A$, cancelling 2 on both sides we get $b(g)=(1-g)a$, so $b\in B^1(G,A)$.
\item The map $H^1(G/\langle c\rangle, A)\rightarrow H^1(G,A)$ is clearly injective. For surjectivity, fix $b\in Z^1(G,A)$ and consider $b(c^2)=b(1)=0$: using the cocycle relation
$$0=(1+c)b(c)=2b(c),$$
hence $b(c)=0$ and $b$ factors through $G/\langle c\rangle$.
\eprsk
\end{enumerate}

\subsection{About $SL_2(\Z)$ and $GL_2(\Z)$}

\medskip
We denote by $C_n$ the cyclic group of order $n$, and by $D_n$ the dihedral group of order $2n$. Set 
$$
s=\left(\begin{array}{cc}0 & -1 \\1 & 0\end{array}\right), \quad 
t=\left(\begin{array}{cc}0 & -1 \\1 & 1\end{array}\right), \quad
\varepsilon = 
\begin{pmatrix}
   -1 & 0 \\
   0 & -1 
\end{pmatrix}, \quad
w=\left(\begin{array}{cc}0 & 1 \\1 & 0\end{array}\right) \quad
.$$
It is well-known (see Example 1.5.3 in \cite{Serre}) that $SL_2(\Z)$ admits a decomposition as an amalgamated product:
$$SL_2(\Z)\simeq C_4\ast_{C_2} C_6.$$
with $C_4=\langle s\rangle, C_6=\langle t\rangle$ and $C_2=\{1,\varepsilon\}$.
It extends to an amalgamated product decomposition of $GL_2(\Z)$ (see Section 6 in \cite{BuTa}):
\begin{equation}\label{amalgamGL2}
GL_2(\Z)\simeq D_4\ast_{D_2} D_6,
\end{equation}
with $D_4=\langle s,w\rangle, D_6=\langle t,w\rangle$ and $D_2=\langle\varepsilon, w\rangle$.


\begin{Lem}\label{amenrad} The amenable radical of $GL_2(\Z)$ is the subgroup $C_2$.
\end{Lem}

\bpr This is a very particular case of a result of Cornulier (Proposition 7 in \cite{Cornu}): let $G=A\ast_C B$ be an amalgamated product such that $[A:C]\geq 2$ and $[B:C]\geq 3$. Then the amenable radical of $G$ is the largest normal subgroup of $C$ which is amenable and normalized both by $A$ and $B$.
\eprsk

Our aim now is to describe the maximal amenable subgroups of $PSL_2(\Z)$. Recall that an element $A\in SL_2(\R), A\neq\pm Id$, is {\it elliptic} if $|Tr(A)|<2$,  {\it parabolic} if $|Tr(A)|=2$, and {\it hyperbolic} if $|Tr(A)|>2$. These concepts clearly descend to $PSL_2(\R)$. For $A=\left(\begin{array}{cc}a & b \\c & d\end{array}\right)\in SL_2(\R)$, we denote by $\left[\begin{array}{cc}a & b \\c & d\end{array}\right]$ the image of $A$ in $PSL_2(\R)$.

\begin{Lem}\label{Maxamen}\begin{enumerate}
\item The maximal amenable subgroups ot $PSL_2(\Z)$ are isomorphic to $C_3$, or $\Z$, or $D_\infty$ (the infinite dihedral group).
\item If $A\in PSL_2(\Z)$ is parabolic, then $A$ is contained in a unique maximal amenable subgroup, isomorphic to $\Z$. 
\item If $A=\left[\begin{array}{cc}a & b \\c & d\end{array}\right]\in PSL_2(\Z)$ is hyperbolic, then $A$ is contained in a unique maximal amenable subgroup, isomorphic either to $\Z$ or $D_\infty$. The second case happens if $A$ is conjugate to a symmetric matrix. If the second case happens, then the integer binary quadratic form $$Q_A(x,y)=: bx^2+(d-a)xy -cy^2$$ represents $-b$ over $\Z$.
\end{enumerate}
\end{Lem}

\bpr\begin{enumerate}
\item Since $PSL_2(\Z)$ is the free product $C_2\ast C_3$, by the Kurosh theorem any subgroup of $PSL_2(\Z)$ is a free product of a free group with conjugates of the free factors. So a subgroup does not contain a free group if and only if it is isomorphic to one of the following: $C_2, C_3, \Z, C_2\ast C_2=D_\infty$; note that these four groups are amenable. It remains to show that any copy of $C_2$ is contained in at least one copy of $D_\infty$; but any element of order 2 in $PSL_2(\Z)$ is conjugate to the image of $s$, and the images of $s$ and $tst^{-1}$ together generate a copy of $D_\infty$.
\item If $A$ is parabolic then, up to conjugacy in $PSL_2(\Z)$,  we may assume that $A=\left[\begin{array}{cc}1 & k \\0 & 1\end{array}\right]$ for some $k\in\Z,\,k\neq 0$. Let $H$ be a maximal amenable subgroup of $PSL_2(\Z)$ containing $A$; by the previous point $H$ is isomorphic either to $\Z$ or to $D_\infty$. In both cases $H$ normalizes the subgroup $\langle A\rangle$.  But by direct computation the normalizer of $\langle A\rangle$ is $\langle \left[\begin{array}{cc}1 & 1 \\0 & 1\end{array}\right]\rangle$, which is therefore the unique maximal amenable subgroup contaning $A$.
\item Let $H$ be a maximal amenable subgroup of $PSL_2(\Z)$ containing $A$. Again, $H$ normalizes $\langle A\rangle$. But $A$ is conjugate in $PSL_2(\R)$ to $\left[\begin{array}{cc}\lambda & 0 \\0 & 1/\lambda\end{array}\right]$ (with $\lambda>1$) whose normalizer in $PSL_2(\R)$ is 
$$\R\rtimes C_2=\{\left[\begin{array}{cc}e^t & 0 \\0 & e^{-t}\end{array}\right]: t\in\R\}\rtimes \langle \left[\begin{array}{cc}0 & -1 \\1 & 0\end{array}\right]\rangle.$$
So $H$ is the intersection of $PSL_2(\Z)$ with a given conjugate of $\R\rtimes C_2$, which proves uniqueness of $H$ (and re-proves that $H$ is isomorphic either to $\Z$ or to $D_\infty$).

For the second statement, we may assume that $A$ is symmetric. Then, as the eigenspaces of $A$ are orthogonal in $\R^2$, we see without computation that $\left[\begin{array}{cc}0 & -1 \\1 & 0\end{array}\right]$ conjugates $A$ to $A^{-1}$, so that $A$ and $\left[\begin{array}{cc}0 & -1 \\1 & 0\end{array}\right]$ together generate a copy of $D_\infty$ (then contained in a maximal one). 

For the final statement, observe that the maximal amenable subgroup containing $A$ is isomorphic to $D_\infty$ if and only if there exists an involution $B\in PSL_2(\Z)$ that conjugates $A$ into $A^{-1}$. Assume such a $B$ exists. Then, denoting by $\overline{A},\overline{B}$ lifts of $A,B$ in $SL_2(\Z)$, we will have $\overline{B}\,\overline{A}\,\overline{B}^{-1}=\pm\overline{A}^{-1}$. Taking the trace of both sides, we see that the minus sign leads to $Tr(\overline{A})=0$, in contradiction with $|Tr(\overline{A})|>2$. So $\overline{B}\,\overline{A}\,\overline{B}^{-1}=\overline{A}^{-1}$ or $\overline{B}\,\overline{A}=\overline{A}^{-1}\,\overline{B}$. Now $\overline{B}$ has order 4 in $SL_2(\Z)$, and therefore $Tr(\overline{B})=0$. So we may write $\overline{B}=\left(\begin{array}{cc}x & y \\z & -x\end{array}\right)$ and we get
$$\left(\begin{array}{cc}x & y \\z & -x\end{array}\right)\cdot \left(\begin{array}{cc}a & b \\c & d\end{array}\right)=\left(\begin{array}{cc}d & -b \\-c & a\end{array}\right)\cdot\left(\begin{array}{cc}x & y \\z & -x\end{array}\right).$$
Equating the (1,1)-coefficients on both sides yields:
$$ax+cy=dx-bz \Longleftrightarrow z=\frac{(d-a)x-cy}{b}.$$
Taking into account the condition $\det(\overline{B})=1$, i.e. $-1=x^2+yz$, and inserting the previous value of $z$:
$$-1=x^2+ y(\frac{(d-a)x-cy}{b})\Longleftrightarrow Q_A(x,y)=-b,$$
which concludes the proof. Note that if $A$ is symmetric, i.e. $b=c$, then $Q_A(0,1)=-b$. \eprsk
\end{enumerate}

\begin{Ex} We claim that the maximal amenable subgroup of $PSL_2(\Z)$ containing $A=\left[\begin{array}{cc}3 & 1 \\2 & 1\end{array}\right]$ is infinite cyclic. To see it, we check that the quadratic form $Q_A(x,y)=x^2 -2xy-2y^2$ does not represent -1. But the equation $x^2 -2xy-2y^2=-1$ is equivalent to $(x-y)^2-3y^2=-1$, which leads to Pell's equation $X^2-3Y^2=-1$; as the fundamental unit $2+\sqrt{3}$ in $\Z[\sqrt{3}]$ has norm 1, this equation has no solution.
\end{Ex}

We now lift the results of Lemma \ref{Maxamen} to $SL_2(\Z)$. Since a maximal amenable subgroup of $SL_2(\Z)$ must clearly contain the center, we see that maximal amenable subgroups of $SL_2(\Z)$ are exactly pullbacks of maximal amenable subgroups of $PSL_2(\Z)$ by the quotient map $SL_2(\Z)\rightarrow PSL_2(\Z)$. Observe that the inverse image of an involution in $PSL_2(\Z)$ has order 4 in $SL_2(\Z)$; consequently the inverse image in $SL_2(\Z)$ of a copy of $D_\infty$ in $PSL_2(\Z)$, will be isomorphic to the semi-direct product $\Z\rtimes C_4$ where a generator of $C_4$ acts on $\Z$ by $n\mapsto -n$. So we get immediately the following, that improves on Proposition 2.13 in \cite{JiSk}:

\begin{Prop}\label{MaxAmen}\begin{enumerate}
\item The maximal amenable subgroups ot $SL_2(\Z)$ are isomorphic to $C_6$, or $\Z\times C_2$, or $\Z\rtimes C_4$.
\item For parabolic (resp. hyperbolic) matrices in $SL_2(\Z)$, item (2) (resp. item (3)) of lemma \ref{Maxamen} applies with the obvious changes. \eprsk
\end{enumerate}
\end{Prop}

\section{Relative property (T)}

\begin{Prop}\label{BurgerT} Let $G$ be a non-amenable subgroup of $GL_2(\Z)$.
\begin{enumerate}
\item The restriction  of the representation $\rho_n$ to $G$ is irreducible.
\item The pair $(P_n(\Z)\rtimes_{\rho_n}G, P_n(\Z))$ has the relative property (T). In particular $P_n(\Z)\rtimes_{\rho_n}G$ is not Haagerup.
\end{enumerate}
\end{Prop}

\bpr \begin{enumerate}
\item Since $G\cap SL_2(\Z)$ has index at most 2 in $G$, replacing $G$ by $G\cap SL_2(\Z)$ we may assume that $G\subset SL_2(\Z)$. Let $L$ be the Zariski closure of $G$ in $SL_2(\R)$, so that $L$ is a Lie subgroup of $SL_2(\R)$, hence of dimension 0, 1, 2 or 3. As Lie subgroups of dimension 0, 1, 2 are virtually solvable hence amenable, $L$ has dimension 3, i.e. $G$ is Zariski dense in $SL_2(\R)$. Since the representation $\rho_n$ is algebraic, irreducibility is preserved by passing to a Zariski dense subgroup.\footnote{We recall the argument: if $W$ is a $\rho_n(G)$-invariant subspace, and $(f_i)_{i\in I}$ is a set of linear forms such that $W=\cap_{i\in I}\Ker(f_i)$, then $\rho_n(G)$-invariance of $W$ is equivalent to $f_i(\rho_n(g)w)=0$ for every $g\in G, w\in W, i\in I$. View this as a system of polynomial equations in the matrix coefficients of $g$: it vanishes on $G$, hence also vanishes on $SL_2(\R)$ by Zariski density.}

\item Set $V_n=P_n(\R)$. By Proposition 7 in \cite{Burger} (see especially Example (2) on p. 62 of \cite{Burger}), if $G$ does not fix any probability measure on the projective space $P(V_n^*)$, then the pair $(P_n(\Z)\rtimes_{\rho_n}G,P_n(\Z))$ has the relative property (T). Since the representation $\rho_n$ of $SL_2(\R)$ is equivalent to its contragredient $\rho_n^*$, it is enough to check that there is no $G$-fixed probability measure on the projective space $P(V_n)$. So assume by contradiction that there is such a probability measure $\mu$. Then, by Corollary 3.2.2 in \cite{Zim}, there are exactly two cases:
\begin{itemize}
\item The measure $\mu$ is not supported on a finite union of proper projective subspaces. Then the stabilizer $PGL(V_n)_\mu$ is compact, which contradicts the fact that the image of $G$ in $PGL(V_n)$ is infinite discrete.
\item There exists a proper linear subspace $W$ in $V_n$ such that $\mu([W])>0$ (where $[W]$ denotes the image of $W$ in $P(V_n)$), and moreover the orbit of $[W]$ under the stabilizer $PGL(V_n)_\mu$ is finite. In particular there is a finite index subgroup of $PGL(V_n)_\mu$ that leaves $[W]$ invariant. So there is a finite index subgroup $G_0$ of $G$ that leaves the linear subspace $W$ invariant, in contradiction with the first part of the Proposition.
\eprsk
\end{itemize}
\end{enumerate}

\section{Maximal Haagerup subgroups}

An interesting question raised in \cite{JiSk} is whether every countable group admits a Haagerup radical, i.e. a unique maximal normal subgroup with the Haagerup property. We first show that $G_n$ admits such a Haagerup radical.

\begin{Prop} For every $n\geq 1$, the Haagerup radical of $G_n$ is $P_n(\Z)\rtimes_{\rho_n} C_2$.
\end{Prop}

\bpr Set $U=:P_n(\Z)\rtimes_{\rho_n}C_2$. Let $N\triangleleft G_n$ be a normal Haagerup subgroup, we want to prove that $N$ is contained in $U$. We proceed as in the proof of Proposition 2.10 in \cite{JiSk}: the subgroup $UN$ is normal and since $UN/N\simeq U/(U\cap N)$ is amenable, $UN$ is an amenable extension of a Haagerup group, hence $UN$ is Haagerup. Since $UN$ contains $U$ we have in particular $P_n(\Z)\subset UN$ so $UN=\Z^{n+1}\rtimes_{\rho_n}K$, for some normal subgroup $K\triangleleft GL_2(\Z)$. By Proposition \ref{BurgerT}(2), the subgroup $K$ must be amenable, i.e. $K\subset C_2$ by lemma \ref{amenrad}. So $UN\subset U$ and therefore $N\subset U$.
\eprsk

Since $P_n(\Z)\rtimes_{\rho_n} C_2$ is actually amenable, we immediately have:

\begin{Cor} The amenable radical of $G_n$ is $P_n(\Z)\rtimes_{\rho_n} C_2$.\epr
\end{Cor}

\medskip
We now come to the proof of our first main result.


\bigskip\noindent
{\bf Proof of Theorem \ref{main}:} Let 
$$q_n: P_n(\Z)\rtimes GL_2(\Z)\rightarrow GL_2(\Z): (v,S)\mapsto S$$
 be the quotient map. Observe that, as $\Ker(q_n|_H)=H\cap P_n(\Z)$ is abelian, $H$ is amenable if and only if $q_n(H)$ is amenable. We separate the two cases.
\begin{enumerate}
\item If $H$ is amenable, set $C=q_n(H)$, so that $H$ is contained in the amenable subgroup $P_n(\Z)\rtimes_{\rho_n}C$. By maximality, we have $H=P_n(\Z)\rtimes_{\rho_n} C$, and $C$ is maximal amenable in $G$.

\item Assume now that $H$ is not amenable, and set $K=q_n(H)$.

{\it Claim 1:} the subgroup $H\cap P_n(\Z)$ is invariant by $\rho_n(K)$. Indeed, fix $h\in H$ and write $h=(v_h,q_n(h))$ as an element of the semi-direct product $G_n$. For $(w,1)\in H\cap P_n(\Z)$, we have, because $P_n(\Z)$ is abelian:
$$(\rho_n(q_n(h^{-1}))w,1)=(0,q_n(h))^{-1}(w,1)(0,q_n(h))$$
$$=(0,q_n(h))^{-1}(-v_h,1)(w,1)(v_h,1)(0,q_n(h))$$
$$=((v_h,1)(0,q_n(h)))^{-1}(w,1)((v_h,1)(0,q_n(h)))$$
$$=(v_h,q_n(h))^{-1}(w,1)(v_h,q_n(h))=h^{-1}(w,1)h$$
which belongs to $H\cap P_n(\Z)$ because the latter is normal in $H$. This proves Claim 1.

{\it Claim 2:} We have $H\cap P_n(\Z)=\{0\}$. To see it, let $k$ be the rank of the free abelian group $H\cap P_n(\Z)$, so that $0\leq k\leq n+1$, we must show that $k=0$. 

If $k=n+1$, then $H\cap P_n(\Z)$ has finite index in $P_n(\Z)$, so that $H$ has finite index in $H\cdot P_n(\Z)$. By maximality, we must have $H\cdot P_n(\Z)=H$, i.e. $P_n(\Z) \subset H$ and $H= P_n(\Z)\rtimes_{\rho_n}K$; as $K$ is not amenable, this contradicts part (2) of Proposition \ref{BurgerT}.

If $1\leq k\leq n$, we denote by $W$ the linear subspace of $P_n(\R)$ generated by $H\cap\Z^{n+1}$. By Claim 1, the subspace $W$ is invariant by $\rho_n(K)$, contradicting Part 1 of Proposition \ref{BurgerT}. This proves Claim 2.

At this point we know that $q_n|_H$ induces an isomorphism from $H$ onto $K$, so that by Proposition \ref{1cohom} there exists a 1-cocycle $b\in Z^1(K,P_n(\Z))$ so that $H=\{(b(k),k):k\in K\}$. By maximality of $H$, the 1-cocycle cannot be extended to a larger subgroup of $G$. This proves the direct implication of the Theorem.

\end{enumerate}

For the converse, if $C$ is a maximal amenable subgroup of $G$, then $P_n(\Z)\rtimes C$ is Haagerup, and maximality follows immediately from Part 2 of Proposition \ref{BurgerT}. If $K$ is a non-amenable subgroup of $G$ and $b\in Z^1(K,P_n(\Z))$ is a 1-cocycle that cannot be extended to a larger subgroup, then $H=\{(b(k),k):k\in K\}$ is a Haagerup subgroup of $P_n(\Z)\rtimes_{\rho_n} G$, and maximality follows the dichotomy in the direct implication of the Theorem.
\eprsk

\section{Cohomological matters: $n$ odd}

\subsection{Maximal subgroups in $PSL_2(\Z)$}

Recall that a subgroup $H$ in a group $G$ is said to be {\bf maximal} if $H$ is maximal among proper subgroups of $G$. 

\begin{Prop}\label{maximalfree} The group $PSL_2(\Z)$ admits free maximal subgroups of arbitrary large finite index; and also of infinite index.
\end{Prop}

\bpr Since $PSL_2(\Z)$ is isomorphic to the free product $C_2\ast C_3$, a subgroup $H$ of $PSL_2(\Z)$ is free if and only it is torsion-free, if and only if it has no element of order 2 or 3, if and only if it does not meet the conjugacy classes of (the images of) $s$ and $t$ in $PSL_2(\Z)$.

\begin{enumerate}
\item {\it Finite index}. Let $p$ be a prime, with $p\equiv 11 \mod 12$. Consider the congruence subgroup
$$\overline{\Gamma}_0(p)=\{\left(\begin{array}{cc}a & b \\c & d\end{array}\right)\in PSL_2(\Z):c\equiv 0 \mod p\}.$$
So $\overline{\Gamma}_0(p)$ is the inverse image under reduction modulo $p$, of the upper triangular subgroup in $PSL_2(p)$. Hence $\overline{\Gamma}_0(p)$ is a maximal subgroup of index $p+1$ in $PSL_2(\Z)$.
For $g=\left(\begin{array}{cc}a & b \\c & d\end{array}\right)\in PSL_2(\Z)$, by direct computation we have that the (2,1)-entry of $gsg^{-1}$ is $c^2+d^2$, which is not divisible by $p$ as $p\equiv 3 \mod 4$; so $\overline{\Gamma}_0(p)$ has no element of order 2. Similarly, the (2,1)-entry of $gtg^{-1}$ is $c^2-cd+d^2$. Assume by contradiction that $p$ divides $c^2-cd+d^2$. Then it also divides $4(c^2-cd+d^2)=(2c-d)^2+3d^2$, from which it follows that the Legendre symbol $\left(\begin{array}{c}-3\\
\hline
 p\end{array}\right)=-1$, as one sees by inverting $d$ modulo $p$. Using quadratic reciprocity, this contradicts $p\equiv 2 \mod 3$. So $\overline{\Gamma}_0(p)$ has no element of order 3, and therefore is a free group. 
 
 \item {\it Infinite index}. We will appeal to Theorem E in \cite{FLMS} together with its proof. Consider the following two permutations of $\N=\{1,2,...\}$:
 $$a=(12)(34)(56)(78)...$$
 $$b=(123)(456)(789)...$$
 View the symmetric group $Sym(\N)$ as a Polish group with the topology of pointwise convergence. Let $T$ be the subset of $Sym(\N)$ consisting of those $\vartheta$'s in $Sym(\N)$ such that the action of $PSL_2(\Z)$ on $\N$ defined by mapping $s$ to $a$ and $t$ to $\vartheta b\vartheta^{-1}$, is transitive on $\N$. By lemma 6.2 in \cite{FLMS}, the set $T$ is a non-empty $G_\delta$ in $Sym(\N)$, hence $T$ is itself Polish. Now Theorem E together with its proof in section 6.2 of \cite{FLMS}, shows that for a generic choice of $\vartheta\in T$, the corresponding action of $PSL_2(\Z)$ will be 2-transitive (in fact highly transitive, i.e. $n$-transitive for every $n\geq 1$). By 2-transitivity, the stabilizer of any point in $\N$ will be a maximal subgroup of infinite index in $PSL_2(\Z)$. Since $a$ and $\vartheta b\vartheta^{-1}$ are permutations without fixed points, that stabilizer avoids the conjugacy classes of $s$ and $t$, and therefore it is free.
 \eprsk
\end{enumerate}

\subsection{The case $n=1$}

The next result is specific to $n=1$; it applies in particular to $G=SL_2(\Z)$ and $G=GL_2(\Z)$. For $SL_2(\Z)$, a different proof appears in lemma 2.14 in \cite{JiSk}.

\begin{Prop}\label{order6} Let $G$ be a subgroup of $GL_2(\Z)$  containing an element $t_0$ of order 6. Let $\rho:G \rightarrow GL_N(\Z)$ be a homomorphism such that $\rho(t_0^3)=-1$, and $\rho(t_0)$ does not admit $-1$ as an eigenvalue. Then $H^1(G,\Z^N)=0$. This applies in particular to $\rho=\rho_1|_G$.
\end{Prop}

\bpr In an amalgamated product of groups, any element of finite order is conjugate to a finite order element in one of the factors (see e.g. Cor. 1 in section 3 of \cite{Serre}). So in $GL_2(\Z)= D_4\ast_{D_2} D_6$, conjugating $G$ if necessary we may assume that $t_0=t$, so that $\varepsilon=t_0^3$ belongs to $G$. 

Fix $b\in Z^1(G,\Z^N)$, we want to prove that $b$ is a 1-coboundary. By lemma \ref{central} (applied with $z=\varepsilon$) it is enough to prove that $b(\varepsilon)\in 2\Z^N$.

 Set $T=:\rho(t_0)$, then:
$$0=T^3+1=(T+1)(T^2-T+1).$$
Since by assumption $T+1$ is invertible in $M_N(\C)$, we have $0=T^2-T+1$, i.e. $2T=T^2+T+1$. Now expanding $t_0^3=\varepsilon$ by the 1-cocycle relation we get
$$b(\varepsilon)=(T^2+T+1)b(t_0)=2Tb(t_0).$$
i.e. $b(\varepsilon)\in 2\Z^N$. For the final statement about $\rho_1$, it is enough to observe that $-1$ is not an eigenvalue of $t=\left(\begin{array}{cc}0 & -1 \\1 & 1\end{array}\right)$.
\eprsk 


We now revisit an example from Jiang-Skalski (Proposition 2.16 in \cite{JiSk}): for $N\geq 2$, define
$$\Gamma_1(N)=\{\left(\begin{array}{cc}a & b \\c & d\end{array}\right)\in SL_2(\Z):c\equiv 0 \mod N,\,a\equiv d\equiv 1\mod N\}.$$
Define also the vector $v_N=\left(\begin{array}{c}1/N \\0\end{array}\right)$ in $\Q^2$. Then Jiang-Skalski prove that $\Gamma_1(N)$ is exactly the set of elements of $SL_2(\Z)$ such that $b_N(g)=gv_N-v_N$ belongs to $\Z^2$, so that $b_N$ is a 1-cocycle in $Z^1(\Gamma_1(N),\Z^2)$ that does not extend to an overgroup, and $H_N=\{(b_N(k),k):k\in \Gamma_1(N)\}$ is maximal Haagerup in $\Z^2\rtimes SL_2(\Z)$. We make the result more precise by observing that $H_N$ is actually a free group for infinitely many values of $N$. 

\begin{Prop} If $N$ is divisible by some prime $p$ with $p\equiv 11\mod 12$, then $H_N$ is a free group.
\end{Prop}

\bpr We will use the fact that any subgroup of a free group is free. It is enough to show that $\Gamma_1(N)$ is free. Since $\Gamma_1(N)\subset\Gamma_1(p)$, we reduce to the case of $\Gamma_1(p)$. Clearly $\varepsilon\notin \Gamma_1(p)$, so $\Gamma_1(p)$ maps injectively to its image in $PSL_2(\Z)$. Clearly the latter is contained in $\overline{\Gamma}_0(p)$, which is free by the proof of Proposition \ref{maximalfree}. Hence the conclusion.
\eprsk

\subsection{$n$ odd, general case}

For general odd $n$, we have the following:

\begin{Prop}\label{finite2-group} Let $G$ be a finitely generated subgroup of $GL_2(\Z)$. Assume that $G$ is generated by $k$ elements $g_1,g_2,..., g_k$ and contains $\varepsilon = 
\begin{pmatrix}
   -1 & 0 \\
   0 & -1 
\end{pmatrix}$. Then, for $n$ odd, $H^1(G,P_n(\Z))$ is a vector space over the field $\F_2$ of two elements, with dimension at most $(k-1)(n+1)$.
\end{Prop}

\bpr By the first part of lemma \ref{central}, we know that $H^1(G,P_n(\Z))$ is a vector space over $\F_2$. Now any cocycle $b\in Z^1(G,P_n(\Z))$ is completely determined by the vectors $b(g_1),..., b(g_k)$. On the other hand the boundary map $\partial: P_n(\Z)\rightarrow Z^1(G,P_n(\Z)): P\mapsto \partial P$, with $(\partial P)(g)=\rho_n(g)P-P$, is injective. So $H^1(G,P_n(\Z))$ appears as a sub-quotient of $(P_n(\Z))^{k-1}$, which gives the desired bound on the dimension.
\eprsk

\medskip
{\bf Proof of Theorem \ref{freesbgp}:} Let $M$ be a a maximal free subgroup of $PSL_2(\Z)$ as in Proposition \ref{maximalfree}. Lifting arbitrarily a free basis of $M$ to $SL_2(\Z)$, we get a free subgroup $K$ in $SL_2(\Z)$, say of rank $k$ (possibly $k=\infty$). In particular, the rank of $H^1(K,P_n(\Z))$ as an abelian group is $(k-1)(n+1)$. Now by maximality of $M$ in $PSL_2(\Z)$, together with the fact that the central extension
$$1\rightarrow\langle\varepsilon\rangle \rightarrow SL_2(\Z)\rightarrow PSL_2(\Z)\rightarrow 1$$
does not split, we see that the only overgroups of $K$ in $SL_2(\Z)$ are $K\times\langle\varepsilon\rangle$ and $SL_2(\Z)$. For both of them, $H^1$ is a torsion group (by Proposition \ref{finite2-group}), so any 1-cocycle $b$ of infinite order in $H^1(K,P_n(\Z))$ will not extend to an overgroup of $K$, and the subgroup $H=\{(b(k),k):k\in K\}$ is maximal Haagerup in $P_n(\Z)\rtimes SL_2(\Z)$.
\eprsk

\section{Cohomological matters: $n$ even}

We will see that for $n$ even the groups $H^1(SL_2(\Z),P_n(\Z))$ and $H^1(GL_2(\Z),P_n(\Z))$ are infinite, so that the situation is completely different from the case $n$ odd.

Since $\rho_n(\varepsilon)=1$ as $n$ is even, by the second part of Lemma \ref{central} we may replace $H^1(SL_2(\Z),P_n(\Z))$ by $H^1(PSL_2(\Z),P_n(\Z))$, and similarly $H^1(GL_2(\Z),P_n(\Z))$ by $H^1(PGL_2(\Z),P_n(\Z))$.

\subsection{The case of $SL_2(\Z)$}

Set $S=\rho_n(s),T=\rho_n(t)$, and recall that $S,T$ generate $\rho_n(SL_2(\Z))\simeq PSL_2(\Z)$.

\begin{Lem}\label{cocextends} An assignment $S\mapsto b(S), T\mapsto b(T)$ of two vectors in $P_n(\Z)$, extends to $b\in Z^1(PSL_2(\Z),P_n(\Z))$ if and only if $(1+S)b(S)=0$ and $(1+T+T^2)b(T)=0$.
\end{Lem}

\bpr View $PSL_2(\Z)$ as the free product $\langle S\rangle\ast\langle T\rangle=\langle S,T|S^2=1,T^3=1\rangle$. So the assignment $S\mapsto b(S), T\mapsto b(T)$ extends to a 1-cocycle if and only it satisfies $(S+1)b(S)=0$ (obtained by expanding $S^2=1$ by the 1-cocycle relation) and $(1+T+T^2)b(T)=0$ (obtained by expanding $T^3=1$ by the cocycle relation).
\eprsk

Define, for $n$ even:
$$\eta(n)= \left\{\begin{array}{ccc}0 & if & n\equiv 2 \mod 3 \\1 & if & n\equiv 0 \mod 3 \\-1 & if & n\equiv 1 \mod 3\end{array}\right.$$


\begin{Thm}\label{CohomSL_2}  For $n$ even: 
$$rk(H^1(PSL_2(\Z),P_n(\Z))=\frac{n+1+3(-1)^{\frac{n}{2}+1}-4\eta(n)}{6}>0.$$
In particular $rk(H^1(PSL_2(\Z),P_n(\Z))\geq\frac{n-6}{6}$.
\end{Thm}

\bpr Since we are only interested in the torsion free part of $H^1(PSL_2(\Z),P_n(\Z))$ we may work rationally and compute the $\Q$-dimension of $H^1(PSL_2(\Z),P_n(\Q))$.

We define a subgroup of $Z^1(PSL_2(\Z),P_n(\Q))$ by 
$$Z_0^1=:\{b\in Z^1(PSL_2(\Z),P_n(\Q)): b(T)=0\},$$
We first observe that every cocycle $b\in H^1(PSL_2(\Z),P_n(\Q))$ is cohomologous to a cocycle in $Z^1_0$. This is clear in terms of affine actions: if $\alpha(g)v=gv+b(g)$ (with $v\in P_n(\Q)$) is the affine action associated with $b$, then for any $w\in P_n(\Q)$ the vector $w_0=:\frac{1}{3}(w+\alpha(T)w+\alpha(T^2)w)$ is $\alpha(T)$-fixed, so the coboundary $b'(g)=b(g)+gw_0-w_0$ vanishes on $T$. 

Consequently, setting $m_0=:\dim Z_0^1$ and $n_0=:\dim(Z_0^1\cap B^1(PSL_2(\Z),P_n(\Q)))$, we have 
$$\dim (H^1(PSL_2(\Z),P_n(\Q)))= \dim(Z_0^1/(Z_0^1\cap B^1(PSL_2(\Z),P_n(\Q))))=m_0-n_0,$$
and it is enough to compute $m_0$ and $n_0$ separately. This will be done in 3 steps. Note that, by lemma \ref{cocextends}, the space $Z^1_0$ can be identified with $\ker(S+1)$

\begin{enumerate}

\item {\it Computation of $m_0=\dim (\ker(S+1))$}. Observe that $(SP)(X,Y)=P(Y,-X)$ for $P\in P_n(\Z)$. Expanding $P$ as a sum of monomials $P(X,Y)=\sum_{k=0}^n a_kX^{n-k}Y^k$, we get 
$$(SP)(X,Y)=\sum_{k=0}^n (-1)^k a_k Y^{n-k}X^k=\sum_{k=0}^n (-1)^k a_{n-k}X^{n-k}Y^k,$$
so that $SP=-P$ if and only if 
\begin{equation}\label{1+S}
a_k=(-1)^{k+1}a_{n-k} \:\mbox{for every} \;k=0,1,...,n. 
\end{equation}
So we may choose $a_0,a_1,..., a_{n/2-1}$ arbitrarily, while $a_{n/2}=0$ if $n\equiv 0 \mod 4$, and $a_{n/2}$ can be chosen arbitrarily if $n\equiv 2 \mod 4$. So
\begin{equation}\label{ker(S+1)}
m_0=\frac{n+1+(-1)^{\frac{n}{2}+1}}{2}.
\end{equation}

\item {\it We claim that $Z_0^1\cap B^1(PSL_2(\Z),P_n(\Q))$ can be identified with $\ker(T-1)$}. Indeed $b\in Z^1_0$ is a 1-coboundary if and only if there exists $P\in P_n(\Q)$ such that $b(S)=(S-1)(P)$ and $0=(T-1)(P)$. This is still equivalent to $b(S)\in (S-1)(\ker(T-1))$. This already identifies $Z_0^1\cap B^1(PSL_2(\Z),P_n(\Q))$ with $(S-1)(\ker(T-1))$. 


To prove the claim it remains to show that $(S-1)|_{\ker(T-1)}$ is injective, i.e. $\ker(S-1)\cap\ker(T-1)=\{0\}$. But, as $S,T$ generate $PSL_2(\Z)$, the space $\ker(S-1)\cap\ker(T-1)$ is exactly the space of $\rho_n(SL_2(\Z))$-fixed vectors in $P_n(\Q)$, which is 0 by the first part of Proposition \ref{BurgerT}.


\item {\it Computation of $n_0=\dim(\ker(T-1))$}. Since $T^3=1$, the operator $T$ defines a representation $\sigma$ of the cyclic group $C_3$, and $n_0=\dim(\ker(T-1))$ is the multiplicity of the trivial representation in this representation. Recall that $C_3$ has 3 irreducible representations, all of dimension 1: the trivial representation $\chi_0$ and the characters $\chi_\pm$ defined by $\chi_\pm(T)= e^{\pm\frac{2\pi i}{3}}$. Now $\sigma$ is equivalent to the direct sum of $n_0$ copies of $\chi_0$ with $n_+$ copies of $\chi_+$ and $n_-$ copies of $\chi_-$. Note that $n_+=n_-$ because $\sigma$ is a real representation. Computing the character of $\sigma$ we get
$$
n+1=Tr(1)=n_0\chi_0(1)+n_+\chi_+(1)+n_-\chi_-(1)=n_0+2n_+;
$$
$$
Tr(T)= n_0\chi_0(T)+n_+\chi_+(T)+n_-\chi_-(T)=n_0-n_+.
$$
Solving this system for $n_0$ we get:
\begin{equation}\label{ker(T-1)}
rk(\ker(T-1))=n_0=\frac{n+1+2Tr(T)}{3}.
\end{equation}

It remains to compute $Tr(T)$ to get the exact value of $n_0$. We observe that $\rho_n$ extends to a representation of $SL_2(\C)$ on the space $P_n(\C)$. Restricting to the compact torus $T=\{a_\vartheta=\left(\begin{array}{cc}e^{i\vartheta} & 0 \\0 & e^{-i\vartheta}\end{array}\right):\vartheta\in\R\}$ we have the classical formula (see e.g. formula (7.26) in \cite{Hall}):
$$Tr(\rho_n(a_\vartheta))=\frac{\sin((n+1)\vartheta)}{\sin(\vartheta)}.$$
Since $t$ has order 6 in $SL_2(\Z)$, it is conjugate to $a_{\pi/3}$ in $SL_2(\C)$, and therefore:
\begin{equation}\label{TraceT}
Tr(T)=Tr(\rho_n(a_{\pi/3}))=\frac{\sin((n+1)\pi/3)}{\sin(\pi/3)}=\eta(n).
\end{equation}


\end{enumerate}
Using equations (\ref{ker(S+1)}), (\ref{ker(T-1)}) and (\ref{TraceT}), we get the desired result. Note that for $n=2,4,6$ we get $m_0>n_0$ from the following table:
\begin{equation}\label{table}\left.\begin{array}{c|c|c}n & n_0 & m_0 \\
\hline
2 & 1 & 2 \\4 & 1 & 2 \\6 & 3 & 4\end{array}\right.
\end{equation}

\eprsk

\begin{Cor} For $n\geq 2$ even, let $M$ be a free maximal subgroup of $PSL_2(\Z)$ as in Proposition \ref{maximalfree}; let $K\simeq M\times\langle\varepsilon\rangle$ denote the inverse image of $M$ in $SL_2(\Z)$. Then the semi-direct product $P_n(\Z)\rtimes SL_2(\Z)$ contains infinitely many maximal Haagerup subgroups $H$ of the form $H=\{(b(g),g):g\in K\}$, where $b\in Z^1(SL_2(\Z),P_n(\Z))$, which moreover are pairwise not conjugate under $P_n(\Z)$.
\end{Cor}

\bpr Say that the rank of $M$ is $k$ (possibly $k=\infty$). By lemma \ref{central}, we have $rk(H^1(K,P_n(\Z))=rk(M,P_n(\Z))=(k-1)(n+1)$ as $M$ is free, So by Theorem \ref{CohomSL_2}:
$$rk(H^1(K,P_n(\Z))\geq n+1 > \frac{n+1+3(-1)^{\frac{n}{2}+1}-4\eta(n)}{6}=rk(H^1(SL_2(\Z),P_n(\Z)).$$
So taking for $b$ cocycles on $K$ whose classes have infinite order in the co-kernel of the restriction map $H^1(SL_2(\Z),P_n(\Z))\rightarrow H^1(K,P_n(\Z))$, we may construct maximal Haagerup subgroups of the desired form.
\eprsk

\subsection{The case of $GL_2(\Z)$}

We recall the notations $S=\rho_n(s), T=\rho_n(t)$, to which we add $W=:\rho_n(w)$. Note that, for $P\in P_n(\Z)$, we have $(WP)(X,Y)=P(Y,X)$; in particular $WP=P$ (resp. $WP=-P$) means that $P$ is symmetric (resp. antisymmetric).



To estimate the rank of $H^1(GL_2(\Z),P_n(\Z))$, we will need the following representation-theoretic lemma.

\begin{Lem}\label{repD3}
The $\Q$-dimension of the space $P_n^0=:\{P\in P_n(\Q): P\,\mbox{symmetric},TP=P\}$ is $\frac{n+4+2\eta(n)}{6}$.
\end{Lem}

\bpr We observe that $P_n^0$ is invariant under $\rho_n|_{D_3}$, which suggests to use representation theory of the dihedral group $D_3=\langle T,W\rangle$. This group has 3 irreducible representations defined over $\Q$: the trivial character $\chi_0$, the non-trivial character $\chi_1$ defined by $\chi_1(T)=1$ and $\chi_1(W)=-1$, and the 2-dimensional irreducible representation on vectors in $\R^3$ whose 3 coordinates sum up to 0. The character table of $D_3$ is:
$$\begin{array}{c|ccc} & e & T & W \\
\hline
\chi_0 & 1 & 1 & 1 \\
\chi_1 & 1 & 1 & -1 \\
\pi & 2 & -1 & 0\end{array}$$
Write $\rho_n|_{D_3}=n_0\chi_0\oplus n_1\chi_1\oplus n_\pi \pi$. From the character table, it follows that the dimension of $P_n^0$ 
is exactly the multiplicity $n_0$ 
of $\chi_0$ 
in $\rho_n|_{D_3}$. To compute $n_0$ , we will compute the character of $\rho_n|_{D_3}$. Denote by $P_n^s$ (resp. $P_n^a$) the subspace of symmetric (resp. antisymmetric) polynomials in $P_n(\Q)$. Then
$$Tr(W)=\dim(P_n^s)-\dim(P_n^a)=(\frac{n}{2}+1)-\frac{n}{2}=1,$$
so that the character of $\rho_n|_{D_3}$ is given by:
$$n+1=Tr(1)=n_0+n_1 +2n_\pi;$$
$$Tr(T)=n_0+n_1-n_\pi;$$
$$1=Tr(W)=n_0-n_1.$$
Solving for $n_0$, and using formula (\ref{TraceT}), gives the desired result.
\eprsk

\begin{Thm}\label{CohomGL2} For even $n$, we have 
$$rk(H^1(GL_2(\Z),P_n(\Z))=\frac{n-5+3(-1)^{\frac{n}{2}+1}-4\eta(n)}{12}.$$
In particular $rk(H^1(GL_2(\Z),P_n(\Z)))\geq \frac{n-12}{12}$.
\end{Thm}

\bpr As for $SL_2(\Z)$, we will work rationally. We will appeal to a part of the Hochschild-Serre long exact sequence in group cohomology, namely: if $G$ is a group, $N\triangleleft G$ a normal subgroup, and $V$ is a $G$-module with $V^G=0$, then the restriction map $H^1(G,V)\rightarrow H^1(N,V)^{G/N}$ is an isomorphism (see section 8.1 in \cite{Guich}).

We apply this to $V=P_n(\Q),\,G=GL_2(\Z)$ and $N=SL_2(\Z)$, so that $G/N=\langle W\rangle$ has order 2. We must therefore work out the $W$-invariants on $H^1(PSL_2(\Z), P_n(\Q))$. This amounts to finding the $W$-invariants in $Z^1(PSL_2(\Z), P_n(\Q)$ and dividing by the $W$-invariants in $B^1(PSL_2(\Z),P_n(\Q))$. Of course we will enjoy the fact that the $W$-invariants in a subspace of $P_n(\Q)$ are just the symmetric polynomials in that subspace. 

By the beginning of the proof of Theorem \ref{CohomSL_2}, every cocycle in $Z^1(SL_2(\Z),P_n(\Q)$ is cohomologous to a cocycle in $Z^1_0$, which identifies with $\ker(S+1)$. By formula (\ref{1+S}), a polynomial $P(X,Y)=\sum_{k=0}^na_k X^{n-k}Y^k$ is in $\ker(S+1)$ if and only if $a_k=(-1)^{k+1}a_{n-k}$ for every $k=0,1,...,n$. 
So this polynomial $P(X,Y)$ is symmetric if and only if $a_k=0$ for even $k$, and $a_k=a_{n-k}$ for odd $k$. Hence we get:
\begin{equation}\label{Z^100}
\dim(Z^1_{0})^W=\frac{n+1+(-1)^{\frac{n}{2}+1}}{4}.
\end{equation}
Now we turn to the computation of the dimension of $(Z^1_{0}\cap B^1(PGL_2(\Z),P_n(\Q)))^W$. As in the third step of the proof of Theorem \ref{CohomSL_2}, using $SW=WS$ we identify this space first with $(S-1)(\ker(T-1)\cap\ker(W-1))$, second with $\ker(T-1)\cap\ker(W-1)$, which is the subgroup $P_n^0$ from lemma \ref{repD3}: its dimension is $\frac{n+4+2\eta(n)}{6}$.

Summarizing: we have by (\ref{Z^100}):
$$\dim(H^1(GL_2(\Z),P_n(\Z))= \dim(Z^1_{0})^W-\dim(Z^1_{0}\cap B^1(PGL_2(\Z),P_n(\Q)))^W$$
$$=\frac{n+1+(-1)^{\frac{n}{2}+1}}{4}-\frac{n+4+2\eta(n)}{6}$$
$$=\frac{n-5+3(-1)^{\frac{n}{2}+1}-4\eta(n)}{12}.$$
It is clear that this quantity is bounded below by $\frac{n-12}{12}$.
\eprsk

With this we reach the proof of Theorem \ref{main2}.

\medskip
{\bf Proof of Theorem \ref{main2}}: It is enough to show that the co-kernel of the restriction map $H^1(GL_2(\Z),P_n(\Z))\rightarrow H^1(SL_2(\Z),P_n(\Z))$ has rank $>0$. This follows by comparing the formulae in Theorems \ref{CohomSL_2} and \ref{CohomGL2}. Indeed we have 
$$\frac{n+1+3(-1)^{\frac{n}{2}+1}-4\eta(n)}{6}-\frac{n-5+3(-1)^{\frac{n}{2}+1}-4\eta(n)}{12}=\frac{n+7+3(-1)^{\frac{n}{2}+1}-4\eta(n)}{12}>0$$
where the inequality follows from $7+3(-1)^{\frac{n}{2}+1}-4\eta(n)\geq 0$.
\hfill$\square$


\begin{Ex} For $n=2,4$, we can give explicit cocycles witnessing the fact that the image of the restriction map $Rest: H^1(GL_2(\Z), P_n(\Z))\rightarrow H^1(SL_2(\Z), P_n(\Z))$ has infinite index. 

Recall from the proof of Theorem \ref{CohomSL_2} that we defined $Z_0^1=:\{b\in Z^1(PSL_2(\Z),P_n(\Z)): b(T)=0\}$ and identified it naturally with $\ker(1+S)$. For $a\in\Z$, define $b_a\in Z^1(SL_2(\Z),P_n(\Z))$ by prescribing:
$$(b_a(S))(X,Y)=a(X^n-Y^n);\; b_a(T)=0;$$
it follows from (\ref{1+S}) that $b_a(S)\in\ker(1+S)$, hence $b_a\in Z_0^1$. We first show ihat, for $a\neq 0$, the cocycle $b_a$ is non-zero in $H^1$.

Assume that $b_a$ is a coboundary. By the 3rd step of the proof of Theorem \ref{CohomSL_2}, the space $Z^1_0\cap B^1(SL_2(\Z),P_n(\Z))$ identifies naturally with $(S-1)(\ker(T-1))$. For $n=2,4$, we have $rk(\ker(T-1))=1$ by table (\ref{table}). For $n=2$, the space $\ker(T-1)$ is generated by $X^2-XY+Y^2$ (a $t$-invariant quadratic form on $\Z^2$), and by direct computation one gets $(S-1)(X^2-XY+Y^2)=2XY$, forcing $a=0$. For $n=4$, we replace $X^2-XY+Y^2$ by $(X^2-XY+Y^2)^2$ and proceed analogously.

Finally, as $b_a$ is anti-symmetric, for $a\neq 0$ its class in $H^1$ is certainly not in $H^1(SL_2(\Z),P_n(\Z))^W$ which is the image of the restriction map $H^1(GL_2(\Z),P_n(\Z))\rightarrow H^1(SL_2(\Z),P_n(\Z))$.

\end{Ex}

Coming back to Theorem \ref{CohomGL2}, we observe that $rk(H^1(GL_2(\Z),P_n(\Z)))=0$ for $n=2, 4, 6, 8, 12$. This means that $H^1(GL_2(\Z),P_n(\Z))$ is a finite group for these values. We will show that it is not zero. 

\begin{Cor} Fix $n\geq 2$ even. Define 
$$m= \left\{\begin{array}{ccc}n/4 & if & n\equiv 0 \mod 4 \\(n+2)/4 & if & n\equiv 2 \mod 4.\end{array}\right.$$
Then $H^1(GL_2(Z),P_n(\Z))$ has at least order $2^m$.
\end{Cor}

\bpr We start with two observations.
\begin{enumerate}
\item A 1-cocycle $b\in Z^1_0$ such that $b(S)$ is symmetric, extends from $PSL_2(\Z)$ to $PGL_2(\Z)$ by the prescription $b(W)=0$. Indeed 
$PGL_2(\Z)\simeq D_2\ast_{C_2} D_3$, with $C_2=\langle W\rangle, D_2=\langle S,W\rangle, D_3=\langle T,W\rangle$. So by the cocycle relation $b$ extends to $D_3$ by $b|_{D_3}=0$, and to $D_2$ by defining $b(SW)=b(S)$, which is indeed equal to $b(WS)=Wb(S)$ as $b(S)$ is symmetric. So $b$ extends to the amalgamated product. 
\item If a cocycle $b\in Z^1(PGL_2(\Z),P_n(\Z))$ with $b(W)=0$ is a coboundary, then there exists $P\in P_n(\Z)$ such that $(W-1)P=0$ (i.e. $P$ is symmetric) and $b(S)=(1-S)P$. But $SP(X,Y)=P(Y,-X)=P(-X,Y)$ (by symmetry), so writing $P(X,Y)=\sum_{k=0}^{\frac{n}{2}-1}a_k(X^{n-k}Y^k+X^kY^{n-k}) + a_{n/2}X^{\frac{n}{2}}Y^{\frac{n}{2}}$ we get
$$((1-S)P)(X,Y)=\sum_{k=0}^{\frac{n}{2}-1}a_k(1-(-1)^k)(X^{n-k}Y^k+X^kY^{n-k}) + a_{n/2}(1-(-1)^{\frac{n}{2}})X^{\frac{n}{2}}Y^{\frac{n}{2}},$$
so that $(1-S)P$ only has even coefficients. 
\end{enumerate}

We now come to the proof itself. For $\varepsilon\in\{0,1\}^m$, define $b_\varepsilon(T)=0$ and
$$b_\varepsilon(S)(X,Y)=\left\{\begin{array}{ccc}
\sum_{k=1}^m \varepsilon_k(X^{n-2k+1}Y^{2k-1}+X^{2k-1}Y^{n-2k+1})& if & n\equiv 0 \mod 4 \\
\sum_{k=1}^{m-1} \varepsilon_k(X^{n-2k+1}Y^{2k-1}+X^{2k-1}Y^{n-2k+1}) + \varepsilon_m X^{\frac{n}{2}}Y^{\frac{n}{2}}& if & n\equiv 2 \mod 4.\end{array}\right.$$
Observe that $(S+1)b_\varepsilon(S)=0$ by (\ref{1+S}), so that by lemma \ref{cocextends} we get a cocycle $b_\varepsilon \in Z^1_{0}$ with $b_\varepsilon(S)$ symmetric, so by the first observation in the proof we may extend it to $PGL_2(\Z)$ by $b_\varepsilon(W)=0$. 


On the other hand, if $\varepsilon_1,\varepsilon_2$ are distinct elements in $\{0,1\}^m$, then $(b_{\varepsilon_1}-b_{\varepsilon_2})(S)$ has at least one odd coefficient, so $b_{\varepsilon_1}-b_{\varepsilon_2}$ cannot be a coboundary by the 2nd observation above.
\eprsk



\subsection{A combinatorial application}
We end the paper with an amusing combinatorial consequence of our proof of Theorem \ref{CohomSL_2}. For $n$ even, the sum
$$\sum_{k=0}^{n/2} (-1)^k\left(\begin{array}{c}n-k \\ k \end{array}\right),$$
which is an alternating sum on a diagonal in Pascal's triangle, can be computed by combinatorial means (see \cite{Quinn}). We give a representation-theoretic derivation.

\begin{Cor} For even $n$:
$$\sum_{k=0}^{n/2} (-1)^k\left(\begin{array}{c}n-k \\ k \end{array}\right)=\eta(n)$$
\end{Cor}

. 
\bpr In view of formula (\ref{TraceT}), it is enough to prove that the LHS is equal to $Tr(T)=Tr(T^{-1})$. For this we observe that $T^{-1}(P)(X,Y)=P(X-Y,X)$ for $P\in P_n(\R)$. So in the canonical basis  $X^n, X^{n-1}Y,...,XY^{n-1}, Y^n$ we have, for $k=0,1,...,n$:
$$(T^{-1})(X^{n-k}Y^k)= X^k(X-Y)^{n-k}=\sum_{\ell=0}^{n-k}(-1)^\ell\left(\begin{array}{c}n-k \\\ell\end{array}\right)X^{n-\ell}Y^\ell . $$
In the last sum, the term $X^{n-k}Y^k$ does not appear if $k>\frac{n}{2}$, and it appears with coefficient $(-1)^k\left(\begin{array}{c}n-k \\ k \end{array}\right)$ if $0\leq k\leq\frac{n}{2}$. In other words, the $k$-th diagonal coefficient of the matrix of $T^{-1}$ in the canonical basis is
$$\left\{\begin{array}{ccc}0 & if & k>\frac{n}{2} \\(-1)^k\left(\begin{array}{c}n-k \\ k \end{array}\right) & if & k\leq \frac{n}{2}\end{array}\right.$$

This yields the trace of $T^{-1}$:
$$Tr(T^{-1})=\sum_{k=0}^{n/2} (-1)^k\left(\begin{array}{c}n-k \\ k \end{array}\right).$$
\epr



Author's address:\\

\noindent
Institut de math\'ematiques\\
Universit\'e de Neuch\^atel\\
11 Rue Emile Argand - Unimail\\
CH-2000 Neuch\^atel - SUISSE\\
alain.valette@unine.ch

\end{document}